# Local asymptotic minimax risk bounds in a locally asymptotically mixture of normal experiments under asymmetric loss

Debasis Bhattacharya[1] and A. K. Basu[2]

*Visva-Bharati University and Calcutta University*

**Abstract:** Local asymptotic minimax risk bounds in a locally asymptotically mixture of normal family of distributions have been investigated under asymmetric loss functions and the asymptotic distribution of the optimal estimator that attains the bound has been obtained.

## 1. Introduction

There are two broad issues in the asymptotic theory of inference: (i) the problem of finding the limiting distributions of various statistics to be used for the purpose of estimation, tests of hypotheses, construction of confidence regions etc., and (ii) problems associated with questions such as: how good are the estimation and testing procedures based on the statistics under consideration and how to define 'optimality', etc. Le Cam [12] observed that the satisfactory answers to the above questions involve the study of the asymptotic behavior of the likelihood ratios. Le Cam [12] introduced the concept of 'Limit Experiment', which states that if one is interested in studying asymptotic properties such as local asymptotic minimaxity and admissibility for a given sequence of experiments, it is enough to prove the result for the limit of the experiment. Then the corresponding limiting result for the sequence of experiments will follow.

One of the many approaches which are used in asymptotic theory to judge the performance of an estimator is to measure the risk of estimation under an appropriate loss function. The idea of comparing estimators by comparing the associated risks was considered by Wald [19, 20]. Later this idea has been discussed by Hájek [8], Ibragimov and Has'minskii [9] and others. The concept of studying asymptotic efficiency based on large deviations has been recommended by Basu [4] and Bahadur [1, 2]. In the above context it is an interesting problem to obtain a lower bound for the risk in a wide class of competing estimators and then find an estimator which attains the bound. Le Cam [11] obtained several basic results concerning asymptotic properties of risk functions for LAN family of distributions. Jeganathan [10], Basawa and Scott [5], and Le Cam and Yang [13] have extended the results of Le Cam for Locally Asymptotically Mixture of Normal (LAMN) experiments. Basu and Bhattacharya [3] further extended the result for Locally Asymptotically Quadratic (LAQ) family of distributions. A symmetric loss structure (for example,

---

[1]Division of Statistics, Institute of Agriculture, Visva-Bharati University, Santiniketan, India, Pin 731236

[2]Department of Statistics, Calcutta University, 35 B. C. Road, Calcutta, India, Pin 700019

squared error loss) has been used to derive the results in the above mentioned references. But there are situations where the loss can be different for equal amounts of over-estimation and under-estimation, e. g., there exists a natural imbalance in the economic results of estimation errors of the same magnitude and of opposite signs. In such cases symmetric losses may not be appropriate. In this context Bhattacharya et al. [7], Levine and Bhattacharya [15], Rojo [16], Zellner [21] and Varian [18] may be referred to. In these works the authors have used an asymmetric loss, known as the LINEX loss function. Let $\Delta = \hat{\theta} - \theta$, $a \neq 0$ and $b > 0$. The LINEX loss is then defined as:

$$(1.1) \quad l(\Delta) = b[\exp(a\Delta) - a\Delta - 1].$$

Other types of asymmetric loss functions that can be found in the literature are as follows:

$$l(\Delta) = \begin{cases} C_1 \Delta, & \text{for } \Delta \geq 0 \\ -C_2 \Delta, & \text{for } \Delta < 0, \end{cases} \quad C_1, C_2 \text{ are constants,}$$

or

$$l(\Delta) = \begin{cases} \lambda w(\theta) L(\Delta), & \text{for } \Delta \geq 0, \quad \text{(over-estimation)} \\ w(\theta) L(\Delta), & \text{for } \Delta < 0, \quad \text{(under-estimation)}, \end{cases}$$

where '$L$' is typically a symmetric loss function, $\lambda$ is an additional loss (in percentage) due to over-estimation, and $w(\theta)$ is a weight function.

The problem of finding the lower bound for the risk with asymmetric loss functions under the assumption of LAN was discussed by Lepskii [14] and Takagi [17]. In the present work we consider an asymmetric loss function and obtain the local asymptotic minimax risk bounds in a LAMN family of distributions.

The paper is organized as follows: Section 2 introduces the preliminaries and the relevant assumptions required to develop the main result. Section 3 is dedicated to the derivation of the main result. Section 4 contains the concluding remarks and directions for future research.

## 2. Preliminaries

Let $X_1, \ldots, X_n$ be $n$ random variables defined on the probability space $(\mathcal{X}, \mathcal{A}, P_\theta)$ and taking values in $(S, \mathcal{S})$, where $S$ is the Borel subset of a Euclidean space and $\mathcal{S}$ is the $\sigma$-field of Borel subsets of $S$. Let the parameter space be $\Theta$, where $\Theta$ is an open subset of $R^1$. It is assumed that the joint probability law of any finite set of such random variables has some known functional form except for the unknown parameter $\theta$ involved in the distribution. Let $\mathcal{A}_n$ be the $\sigma$-field generated by $X_1, \ldots, X_n$ and let $P_{\theta,n}$ be the restriction of $P_\theta$ to $\mathcal{A}_n$. Let $\theta_0$ be the true value of $\theta$ and let $\theta_n = \theta_0 + \delta_n h$ ($h \in R^1$), where $\delta_n \to 0$ as $n \to \infty$. The sequence $\delta_n$ may depend on $\theta$ but is independent of the observations. It is further assumed that, for each $n \geq 1$, the probability measures $P_{\theta_o,n}$ and $P_{\theta_n,n}$ are mutually absolutely continuous for all $\theta_0$ and $\theta_n$. Then the sequence of likelihood ratios is defined as

$$L_n(\mathbf{X}_n; \theta_0, \theta_n) = L_n(\theta_0, \theta_n) = \frac{dP_{\theta_n,n}}{dP_{\theta_0,n}},$$

where $\mathbf{X}_n = (X_1, \ldots, X_n)$ and the corresponding log-likelihood ratios are defined as

$$\Lambda_n(\theta_0, \theta_n) = \log L_n(\theta_0, \theta_n) = \log \frac{dP_{\theta_n,n}}{dP_{\theta_0,n}}.$$

Throughout the paper the following notation is used: $\phi_y(\mu, \sigma^2)$ represents the normal density with mean $\mu$ and variance $\sigma^2$; the symbol '$\Rightarrow$' denotes convergence in distribution, and the symbol '$\to$' denotes convergence in $P_{\theta_0,n}$ probability.

Now let the sequence of statistical experiments $E_n = \{\mathcal{X}_n, \mathcal{A}_n, P_{\theta,n}\}_{n \geq 1}$ be a locally asymptotically mixture of normals (LAMN) at $\theta_0 \in \Theta$. For the definition of a LAMN experiment the reader is referred to Bhattacharya and Roussas [6]. Then there exist random variables $Z_n$ and $W_n$ ($W_n > 0$ a.s.) such that

$$\Lambda_n(\theta_0, \theta_n) = \log \frac{dP_{\theta_0 + \delta_n h, n}}{dP_{\theta_0, n}} - hZ_n + \frac{1}{2}h^2 W_n \to 0, \tag{2.1}$$

and

$$(Z_n, W_n) \Rightarrow (Z, W) \text{ under } P_{\theta_0, n}, \tag{2.2}$$

where $Z = W^{1/2} G$, $G$ and $W$ are independently distributed, $W > 0$ a.s. and $G \sim N(0, 1)$. Moreover, the distribution of $W$ does not depend on the parameter $h$ (Le Cam and Yang [13]).

The following examples illustrate the different quantities appearing in equations (2.1) and (2.2) and in the subsequent derivations.

**Example 2.1** (An explosive autoregressive process of first order). Let the random variables $X_j, j = 1, 2, \ldots$ satisfy a first order autoregressive model defined by

$$X_j = \theta X_{j-1} + \epsilon_j, X_0 = 0, |\theta| > 1, \tag{2.3}$$

where $\epsilon_j$'s are i.i.d. $N(0, 1)$ random variables. We consider the explosive case where $|\theta| > 1$. For this model we can write

$$f_j(\theta) = f(x_j | x_1, \ldots, x_{j-1}; \theta) \propto e^{-\frac{1}{2}(x_j - \theta x_{j-1})^2}.$$

Let $\theta_0$ be the true value of $\theta$. It can be shown that for the model described in (2.3) we can select the sequence of norming constants $\delta_n = \frac{(\theta_0^2 - 1)}{\theta_0^n}$ so that (2.1) and (2.2) hold. Clearly $\delta_n \to 0$ as $n \to \infty$. We can also obtain $W_n(\theta_0), Z_n(\theta_0)$ and their asymptotic distributions, as $n \to \infty$, as follows:

$$W_n(\theta_0) = \frac{(\theta_0^2 - 1)^2}{\theta_0^{2n}} \sum_{j=1}^n X_{j-1}^2 \Rightarrow W \text{ as } n \to \infty, \text{ where } W \sim \chi_1^2 \text{ and}$$

$$G_n(\theta_0) = \Big(\sum_{j=1}^n X_{j-1}^2\Big)^{-\frac{1}{2}} \Big(\sum_{j=1}^n X_{j-1}\epsilon_j\Big) = \Big(\sum_{j=1}^n X_{j-1}^2\Big)^{\frac{1}{2}} (\hat{\theta}_n - \theta) \Rightarrow G,$$

where $G \sim N(0, 1)$ and $\hat{\theta}_n$ is the m.l.e. of $\theta$. Also

$$Z_n(\theta_0) = W_n^{\frac{1}{2}}(\theta_0) G_n(\theta_0) = \frac{(\theta_0^2 - 1)}{\theta_0^n} \Big(\sum_{j=1}^n X_{j-1}\epsilon_j\Big) \Rightarrow W^{\frac{1}{2}} G = Z,$$

where $W$ is independent of $G$. It also holds that

$$(Z_n(\theta_0), W_n(\theta_0)) \Rightarrow (Z, W).$$

Hence $Z|W \sim N(0, W)$. In general $Z$ is a mixture of normal distributions with $W$ as the mixing variable.

**Example 2.2** (A super-critical Galton–Watson branching process)**.** Let $\{X_0 = 1, X_1, \ldots, X_n\}$ denote successive generation sizes in a super-critical Galton–Watson process with geometric offspring distribution given by

(2.4) $\qquad P(X_1 = j) = \theta^{-1}(1 - \theta^{-1})^{j-1}, \ j = 1, 2, \ldots, 1 < \theta < \infty.$

Here $E(X_1) = \theta$ and $V(X_1) = \sigma^2(\theta) = \theta(\theta - 1)$. For this model we can write

$$f_j(\theta) = f(x_j | x_1, \ldots, x_{j-1}; \theta) = (1 - \frac{1}{\theta})^{x_j - x_{j-1}} (\frac{1}{\theta})^{x_{j-1}}.$$

Let $\theta_0$ be the true value of $\theta$. Here $\delta_n$ can be chosen as $\frac{\sqrt{\theta_0}(\theta_0 - 1)}{\theta_0^{n/2}}$. For this model the random variables $W_n(\theta_0), Z_n(\theta_0)$ and their asymptotic distributions are:

$$W_n(\theta_0) = \frac{(\theta_0 - 1)}{\theta_0^n} \sum_{j=1}^n X_{j-1} \Rightarrow W \text{ as } n \to \infty,$$

where $W$ is an exponential random variable with unit mean. Here

$$G_n(\theta_0) = [\theta_0(\theta_0 - 1)]^{-\frac{1}{2}} (\sum_{j=1}^n X_{j-1})^{-\frac{1}{2}} \sum_{j=1}^n (X_j - \theta_0 X_{j-1}) \Rightarrow G,$$

where $G \sim N(0, 1)$, and for $W$ independent of $G$,

$$Z_n(\theta_0) = W_n^{\frac{1}{2}}(\theta_0) G_n(\theta_0) \Rightarrow W^{\frac{1}{2}} G.$$

It also holds that

$$(Z_n(\theta_0), W_n(\theta_0)) \Rightarrow (W^{\frac{1}{2}} G, W).$$

The decision problem considered here is the risk in the estimation of a parameter $\theta \varepsilon R^1$ using an asymmetric loss function $l(.)$. Throughout the rest of the manuscript the following assumptions apply:

A1 $l(z) \geq 0$ for all $z, z = \hat{\theta} - \theta$
A2 $l(z)$ is non increasing for $z < 0$, non-decreasing for $z > 0$ and $l(0) = 0$.
A3 $\int_{-\infty}^\infty \int_0^\infty l(w^{-\frac{1}{2}} z) e^{-\frac{1}{2} c w z^2} g(w) dw dz < \infty$ for any $c > 0$, where $g(w)$ is the p.d.f. of the random variable $W$.
A4 $\int_{-\infty}^\infty \int_0^\infty w^{\frac{1}{2}} z^2 l(w^{-\frac{1}{2}} d - z) e^{-\frac{1}{2} c w z^2} g(w) dw dz < \infty$, for any $c, d > 0$.
Define $l_a(y) = \min(l(y), a)$, for $0 < a \leq \infty$. This truncated loss makes $l(y)$ bounded if it is not so.
A5 For given $W = w > 0$, $h(\beta, w) = \int_{-\infty}^\infty l(w^{-\frac{1}{2}} \beta - y) \phi_y(0, w^{-1}) dy$ attains its minimum at a unique $\beta = \beta_0(w)$, and $E\beta_0(W))$ is finite.
A6 For given $W = w > 0$, any large $a, b > 0$ and any small $\lambda > 0$ the function $\tilde{h}(\beta, w) = \int_{-\sqrt{b}}^{\sqrt{b}} l_a(w^{-\frac{1}{2}} \beta - y) \phi_y(0, ((1 + \lambda)w)^{-1}) dy$
attains its minimum at $\tilde{\beta}(w) = \tilde{\beta}(a, b, \lambda, w)$, and $E\tilde{\beta}(a, b, \lambda, W) < \infty$.
A7 $\lim_{a \to \infty, b \to \infty, \lambda \to 0} \tilde{\beta}(a, b, \lambda, w) = \beta_0(w)$.
A8 $E(W^{-\frac{1}{2}}) < \infty$.

Note the following:

1. Assumptions A3 and A4 are general assumptions made to ensure the finiteness of the expected loss and other functions. Assumptions A5 and A6 are satisfied, for example, by convex loss functions.

2. If $l(.)$ is symmetric, then $\beta_0(w) = 0 = \tilde{\beta}(a, b, \lambda, w)$.
3. If $l(.)$ is unbounded, then the assumption A8 is replaced by A8' as $E(W^{-1} \times Z^2 l(W^{-\frac{1}{2}}Z)) < \infty$.

Here we will consider a randomized estimator $\xi(Z, W)$ which can be written as $\xi(Z, W, U)$, where $U$ is uniformly distributed on $[0, 1]$ and independent of $Z$ and $W$. The introduction of randomized estimators is justified since the loss function $l(.)$ may not be convex.

## 3. Main result

Under the set of assumptions and notations stated in Section 2 we can have the following generalization of Hájek's result for LAMN experiments under an asymmetric loss structure.

**Lemma 3.1.** *Let $l(.)$ satisfy assumptions* A1 – A7 *and let $Z|W$ be a normal random variable with mean $\theta\sqrt{w} + \beta_0(w)$ and variance 1. Further let $W$ be a random variable with p.d.f. $g(w)$, then for any $\epsilon > 0$ there is an $\alpha = \alpha(\epsilon) > 0$ and a prior density $\pi(\theta)$ so that for any estimator $\xi(Z, W, U)$ satisfying*

$$(3.1) \qquad P_{\theta=0}(|\xi(Z,W,U) - W^{-\frac{1}{2}}Z| > \epsilon) > \epsilon$$

*the Bayes risk $R(\pi, \xi)$ is*

$$
\begin{aligned}
R(\pi, \xi) &= \int \pi(\theta) R(\theta, \xi) d\theta \\
&= \int \pi(\theta) E(l_a(\xi(Z,W,U) - \theta)|\theta) d\theta \\
&\geq \int l(w^{-\frac{1}{2}}\beta_0(w) - y)\phi_y(0, w^{-1})g(w) dy dw + \alpha.
\end{aligned}
\tag{3.2}
$$

*Proof.* Let the prior distribution of $\theta$ be given by $\pi(\theta) = \pi_\sigma(\theta) = (2\pi)^{-\frac{1}{2}}\sigma^{-1}e^{-\frac{\theta^2}{2\sigma^2}}$, $\sigma > 0$, where the variance $\sigma^2$, which depends on $\epsilon$ as defined in (3.1), will be appropriately chosen later. As $\sigma^2 \longrightarrow \infty$, the prior distribution becomes diffuse. The joint distribution of $Z$, $W$ and $\theta$ is given by

$$(3.3) \qquad f(z|w)g(w)\pi(\theta) = (2\pi)^{-1}\sigma^{-1}e^{-\frac{1}{2}(z-(\theta w^{\frac{1}{2}}+\beta_0(w)))^2 - \frac{1}{2}\frac{\theta^2}{\sigma^2}}g(w).$$

The posterior distribution of $\theta$ given $(W, Z)$ is given by $\psi(\theta|w, z)$, where $\psi(\theta|w, z)$ is $N(\frac{w^{\frac{1}{2}}(z-\beta_0(w))}{r(w,\sigma)}, \frac{1}{r(w,\sigma)})$ and the marginal joint distribution of $(Z, W)$ is given by

$$(3.4) \qquad f(z, w) = \phi_z(\beta_0(w), \sigma^2 r(w, \sigma))g(w),$$

where the function $r(s, t) = s + 1/t^2$. Note that the Bayes' estimator of $\theta$ is $\frac{W^{\frac{1}{2}}(Z-\beta_0(W))}{r(W,\sigma)}$ and when the prior distribution is sufficiently diffused, the Bayes' estimator becomes $W^{-\frac{1}{2}}(Z - \beta_0(W))$.

Now let $\epsilon > 0$ be given and consider the following events:

$$|\frac{W^{\frac{1}{2}}(Z-\beta_0(W))}{r(W,\sigma)}| \leq b - \sqrt{b}, \quad |\xi(Z,W,U) - W^{-\frac{1}{2}}Z| > \epsilon,$$

$$|W^{-\frac{1}{2}}(Z - \beta_0(W))| \leq M, \quad \frac{1}{m} = \frac{1}{\sigma^2}(\frac{2M}{\epsilon} - 1) \leq W \leq m.$$

Then

$$
(3.5) \quad |W^{-\frac{1}{2}}(Z - \beta_0(W)) - \frac{W^{\frac{1}{2}}(Z - \beta_0(W))}{r(W,\sigma)}| = |\frac{\frac{W^{-\frac{1}{2}}(Z-\beta_0(W))}{\sigma^2}}{r(W,\sigma)}|
$$
$$
\leq \frac{M}{\sigma^2 r(W,\sigma)} = \frac{M}{\sigma^2 W + 1} \leq \frac{\epsilon}{2}.
$$

Now, for any large $a, b > 0$, we have

$$
(3.6) \quad \int_{-b}^{b} l_a(\xi(z,w,u) - \theta)\psi(\theta|z,w)d\theta
$$
$$
= \int_{-b}^{b} l_a(\xi(z,w,u) - y - \frac{w^{\frac{1}{2}}(z - \beta_0(w))}{r(w,\sigma)})\phi_y(0, \frac{1}{r(w,\sigma)})dy,
$$

where $y = \theta - \frac{w^{\frac{1}{2}}(z-\beta_0(w))}{r(w,\sigma)}$. Now, since $\theta|z, w \sim N(\frac{w^{\frac{1}{2}}(z-\beta_0(w))}{r(w,\sigma)}, \frac{1}{r(w,\sigma)})$, we have $y|z, w \sim N(0, \frac{1}{r(w,\sigma)})$. It can be seen that $|\xi(z,w,u) - \frac{w^{\frac{1}{2}}(z-\beta_0(w))}{r(w,\sigma)} - w^{-\frac{1}{2}}\beta_0(w)| > \frac{\epsilon}{2}$. Hence due to the nature of the loss function, for a given $w > 0$, we can have, from (3.6),

$$
(3.7) \quad \int_{-b}^{b} l_a(\xi(z,w,u) - \frac{w^{\frac{1}{2}}(z - \beta_0(w))}{r(w,\sigma)} - y)\phi_y(0, \frac{1}{r(w,\sigma)})dy
$$
$$
\geq \int_{-\sqrt{b}}^{\sqrt{b}} l_a(w^{-\frac{1}{2}}\beta_0(w) - y)\phi_y(0, \frac{1}{r(w,\sigma)})dy
$$
$$
\geq \int_{-\sqrt{b}}^{\sqrt{b}} l_a(w^{-\frac{1}{2}}\tilde{\beta}(a,b,\lambda,w) - y)\phi_y(0, \frac{1}{r(w,\sigma)})dy + \delta
$$
$$
= \tilde{h}(\tilde{\beta}(a,b,\lambda,w)) + \delta,
$$

where $\delta > 0$ depends only on $\epsilon$ but not on $a, b, \sigma^2$ and (3.7) holds for sufficiently large $a, b, \sigma^2$ (here $\lambda = \frac{1}{w\sigma^2} \to 0$ as $\sigma^2 \to \infty$ and $\frac{1}{m} \leq w \leq m$).

A simple calculation yields

$$
(3.8) \quad \tilde{h}(\tilde{\beta}(a,b,\lambda,w))
$$
$$
= \int_{-\sqrt{b}}^{\sqrt{b}} l_a(w^{-\frac{1}{2}}\tilde{\beta}(a,b,\lambda,w) - y)\phi_y(0, \frac{1}{r(w,\sigma)})dy
$$
$$
\geq \int_{-\sqrt{b}}^{\sqrt{b}} l_a(w^{-\frac{1}{2}}\tilde{\beta}(a,b,\lambda,w) - y)\phi_y(0, \frac{1}{w})(1 - \frac{y^2}{\sigma^2})dy
$$
$$
= h(\beta_0(w)) - \int_{-\sqrt{b}}^{\sqrt{b}} l_a(w^{-\frac{1}{2}}\tilde{\beta}(a,b,\lambda,w) - y)\frac{y^2}{\sigma^2}\phi_y(0, \frac{1}{w})dy.
$$

Hence

$$R(\pi(\theta), \xi) = \int_{-\infty}^{\infty} \pi(\theta) R(\theta, \xi) d\theta$$

$$\geq \int_{-b}^{b} \pi(\theta) E(l_a(\xi(Z, W, U) - \theta)) d\theta$$

(3.9)
$$= \int_{\theta=-b}^{b} \int_{u=0}^{1} \int_{w=0}^{\infty} \int_{z=-\infty}^{\infty} l_a(\xi(z, w, u) - \theta) \psi(\theta|z, w) f(z, w) d\theta du dw dz$$

$$\geq \int h(\beta_0(w)) g(w) dw \times P(|W^{-\frac{1}{2}}(Z - \beta_0(W))| \geq b - \sqrt{b}) - \frac{k}{\sigma^2}$$

$$+ \delta P\{|\xi(Z, W, U) - W^{-\frac{1}{2}} Z| > \epsilon, |W^{-\frac{1}{2}}(Z - \beta_0(W))| \leq M, \frac{1}{m} \leq W \leq m\},$$

using (3.7), (3.8) and assumption A4, where $k > 0$ does not depend on $a, b, \sigma^2$. Let

$$A = \{(z, w, u) \in (-\infty, \infty) \times (0, \infty) \times (0, 1) : |\xi(z, w, u) - w^{-\frac{1}{2}} z| > \epsilon,$$

$$|w^{-\frac{1}{2}}(z - \beta_0(w))| \leq M, \frac{1}{m} \leq w \leq m\}.$$

Then $P(A|\theta = 0) > \frac{\epsilon}{2}$ for sufficiently large $M$ due to (3.1). Now under $\theta = 0$ the joint density of $Z$ and $W$ is $\phi_z(\beta_0(w), 1) g(w)$. The overall joint density of Z and W is given in (3.4). The likelihood ratio of the two densities is given by

$$\frac{f(z, w)}{f(z, w|\theta = 0)} = \sigma^{-1} r(w, \sigma)^{-\frac{1}{2}} e^{\frac{1}{2}(z - \beta_0(w))^2 \frac{w}{r(w, \sigma)}}$$

and the ratio is bounded below on $\{(z, w) : |w^{-\frac{1}{2}}(z - \beta_0(w))| \leq M, \frac{1}{m} \leq w \leq m\}$ by $\sigma^{-1} r(m, \sigma)^{-\frac{1}{2}} = \frac{1}{(m\sigma^2 + 1)^{\frac{1}{2}}}$. Finally we have

(3.10)
$$P(A) = \int_{u=0}^{1} \int_A f(z, w, u) dz dw du \geq \frac{1}{(m\sigma^2 + 1)^{\frac{1}{2}}} \frac{\epsilon}{2}.$$

Hence for sufficiently large $m$ and $M$, from (3.9), we have

$$R(\pi(\theta), \xi) \geq \int h(\beta_0(w)) g(w) dw [1 - \frac{\alpha}{2h(\beta_0(w))}] - \frac{k}{\sigma^2} + \delta \frac{\epsilon}{2} \frac{1}{(m\sigma^2 + 1)^{\frac{1}{2}}},$$

assuming $P[|W^{-\frac{1}{2}}(Z - \beta_0(W))| \leq b - \sqrt{b}] \geq 1 - \frac{\alpha}{2h(\beta_0(w))}$. That is,

$$R(\pi(\theta), \xi) \geq \int h(\beta_0(w)) g(w) dw - \frac{\alpha}{2} - \frac{k}{\sigma^2} + \delta \frac{\epsilon}{2} \frac{1}{(m\sigma^2 + 1)^{\frac{1}{2}}}.$$

Putting $\delta \frac{\epsilon}{2} (m\sigma^2 + 1)^{-1/2} - \frac{k}{\sigma^2} = \frac{3\alpha}{2}$ we find $R(\pi(\theta), \xi) \geq \int h(\beta_0(w)) g(w) dw + \alpha$. Hence the proof of the result is complete. □

**Theorem 3.1.** *Suppose that the sequence of experiments $\{E_n\}$ satisfies LAMN conditions at $\theta \in \Theta$ and the loss function $l(.)$ meets the assumptions A1–A8 stated in Section 2. Then for any sequence of estimators $\{T_n\}$ of $\theta$ based on $X_1, \ldots, X_n$ the lower bound of the risk of $\{T_n\}$ is given by*

$$\lim_{\delta \to 0} \liminf_{n \to \infty} \sup_{|\theta - t| < \delta} E_\theta \{l(\delta_n^{-1}(T_n - \theta))\} \geq \int l(\beta_0(w) - y) \phi_y(0, w^{-1}) g(w) dy dw$$

*Furthermore, if the lower bound is attained, then*

$$\delta_n^{-1}(T_n - \theta)) - \frac{W_n^{\frac{1}{2}}(Z_n - \beta_0(W))}{r(W_n, \sigma)} \to 0$$

*or, as* $\sigma^2 \to \infty$

$$\delta_n^{-1}(T_n - \theta) - W_n^{-\frac{1}{2}}(Z_n - \beta_0(W)) \to 0.$$

*Proof.* Since the upper bound of values of a function over a set is at least its mean value on that set, we may write, for sufficiently large n,

$$\sup_{|\theta - t| < \delta} E_\theta\{l(\delta_n^{-1}(T_n - \theta))\} \geq \int_{-b}^{b} \pi(h) E\{l_a(\delta_n^{-1}(T_n - t) - h) | \theta = t + \delta_n h\} dh,$$

whatever the values of constants $a, b$ and the prior density $\pi(h)$ may be. Now, let $Z|w, \theta = t + \delta_n h \sim N(\theta w^{\frac{1}{2}} + \beta_0(w), 1)$. Then we can fix some $\delta > 0$ and choose $a, b$ and $\pi(.)$ in such a way that

$$\int_{-b}^{b} \pi(h) E\{l_a(\xi(Z, W, U) - h) | t + \delta_n h\} dh$$
$$\geq \int l(\beta_0(w) - y) \phi_y(0, w^{-1}) g(w) dy dw - \delta,$$

for any estimator $\xi(Z, W, U)$.

Next we use Lemmas 3.3 and 3.4 of Takagi [17], where we set

$$S_n = \delta_n^{-1}(T_n - t), \quad \Delta_{n,t} = W_n^{-\frac{1}{2}}(Z_n - \beta_0(W)), \text{ and}$$
$$S_n(\Delta_{n,t} = x, U = u) = \inf\{y : P(S_n \leq y | \Delta_{n,t} = x) \geq u\}.$$

Let $F_{n,h}$ = distribution of $S_n$ under $P_{n,h}$, $F_{n,h}^*$ = distribution of $S_n(\Delta_{n,t}, U) = \xi_n(Z_n, W, U)$ under $P_{n,h}$, where $U \sim$ Uniform $(0, 1)$ and is independent of $\Delta_{n,t}$; $G_{n,h}$ is the distribution of $\Delta_{n,t}$ and $G_{n,h}^*$ is the distribution of $\Delta_t = W^{\frac{1}{2}}(Z - \beta_0(W))$. As a consequence of this we have (Takagi [17], p.44)

$$\lim_{n \to \infty} ||F_{n,h} - F_{n,h}^*|| = 0 \text{ and } \lim_{n \to \infty} ||G_{n,h} - G_{n,h}^*|| = 0.$$

Now for any estimator $\xi_n(Z_n, W_n, U) = S_n(\Delta_{n,t}, U)$ and for every $h \varepsilon R^1$ we have

$$|E[l_a(\delta_n^{-1}(T_n - t) - h) | t + \delta_n h] - E[l_a(S_n(\Delta_{n,t}, U) - h) | t + \delta_n h]| \longrightarrow 0$$

and

$$|E[l_a(S_n(\Delta_{n,t}, U) - h) | t + \delta_n h] - E[l_a(\xi_n(Z, W, U) - h) | t + \delta_n h]| \longrightarrow 0.$$

Finally

$$\int_{-b}^{b} \pi(h) E\{l(\delta_n^{-1}(T_n - t) - h) | \theta = t + \delta_n h\} dh$$
$$\geq \int_{-b}^{b} \pi(h) E\{l_a(\xi_n(Z, W, U) - h) | t + \delta_n h\} dh$$
$$\geq \int l(\beta_0(w) - y) \phi_y(0, w^{-1}) g(w) dy dw, \text{ for } n \geq n(a, b, \delta_n, \pi)$$

which proves the result. □

**Example 3.1.** Consider the LINEX loss function as defined by (1.1). It can be seen that $l(\triangle)$ satisfies all the assumptions A1–A7 stated in Section 2. Here a simple calculation will yield

$$h(\beta, w) = b(e^{aw^{-\frac{1}{2}}(\beta + \frac{1}{2}\frac{a}{w^{1/2}})} - aw^{-\frac{1}{2}}\beta - 1),$$

and $h(\beta, w)$ attains its minimum at $\beta_0(w) = -\frac{1}{2}\frac{a}{w^{1/2}}$ and $h(\beta_0, w) = -\frac{b}{2}\frac{a^2}{w}$.

## 4. Concluding remarks

From the results discussed in Le Cam and Yang [13] and Jeganathan [10] it is clear that under symmetric loss structure the results derived in Theorem 3.1 hold with respect to the estimator $W_n^{-\frac{1}{2}}(\theta_0)Z_n(\theta_0)$ and its asymptotic counterpart $W^{-\frac{1}{2}}Z$. Here due to the presence of asymmetry in the loss structure the results derived in Theorem 3.1 hold with respect to the estimator $W_n^{-\frac{1}{2}}(\theta_0)(Z_n(\theta_0) - \beta_0(W)) + \beta_0(W)$ and $W^{-\frac{1}{2}}(Z - \beta_0(W)) + \beta_0(W)$.

Now $W_n^{-\frac{1}{2}}(\theta_0)(Z_n(\theta_0) - \beta_0(W)) \Rightarrow W^{-\frac{1}{2}}(Z - \beta_0(W))$. Hence the asymptotic bias of the estimator under asymmetric loss would be $E(W^{-\frac{1}{2}}(\theta_0)(Z - \beta_0(W)) + \beta_0(W) - \theta) = E(\theta + \beta_0(W) - \theta) = E(\beta_0(W))$.

Consider the model described in Example 2.1. Under the LINEX loss we have $\beta_0(w) = -\frac{a}{2}\frac{1}{w^{1/2}}$ (vide Example 3.1). Here the asymptotic bias of the estimator would be $E(\beta_0(W)) = -\frac{a}{2}E(W^{-\frac{1}{2}})$, which is finite due to Assumption A8.

The results obtained in this paper can be extended in the following two directions: (1) To investigate the case when the experiment is Locally Asymptotically Quadratic (LAQ), and (2) To find the asymptotic minimax lower bound for a sequential estimation scheme under the conditions of LAN, LAMN and LAQ considering asymmetric loss function.

**Acknowledgments.** The authors are indebted to the referees, whose comments and suggestions led to a significant improvement of the paper. The first author is also grateful to the Editor for his support in publishing the article.